\documentclass[aop,MSNbibl,seceqn,dvips]{arximspdf}
\usepackage{graphicx}
\doi{10.1214/13-AOP865} %kopijuoti is PTS
\volume{42}
\issue{5}
\pubyear{2014}
\firstpage{2197}
\lastpage{2206}

\makeatletter
\newtheorem{theorem}{Theorem}
\newproclaim{remark}{Remark}
\newcommand{\eqref}[1]{(\ref{#1})}

\def\EE{\mathsf E}
\def\PP{\mathsf P}
\def\R{\mathbb{R}}

\makeatother

\begin{document}
\begin{frontmatter}

\title{A probabilistic solution to the Stroock--Williams~equation}
\runtitle{\hspace*{-10pt}A probabilistic solution to the Stroock--Williams equation}

\begin{aug}
\author[A]{\fnms{Goran} \snm{Peskir}\corref{}\ead[label=e1]{goran@maths.man.ac.uk}}
\runauthor{G. Peskir}
\affiliation{The University of Manchester}
\address[A]{School of Mathematics\\
The University of Manchester\\
Oxford Road\\
Manchester M13 9PL\\
United Kingdom\\
\printead{e1}} %adresu isvedimo komanda gale!
\end{aug}

% HISTORY:
\received{\smonth{8} \syear{2012}}
\revised{\smonth{6} \syear{2013}}

% ABSTRACT
%
\begin{abstract}
We consider the initial boundary value problem
\begin{eqnarray*}
u_t &=& \mu u_x + \tfrac{1}{2}
u_{xx}\qquad (t>0, x \ge0),
\\
u(0,x) &= &f(x) \qquad (x \ge0),
\\
u_t(t,0) &= & \nu u_x(t,0)\qquad (t>0)
\end{eqnarray*}
of Stroock and Williams [\textit{Comm. Pure Appl. Math.} \textbf{58}
(2005) 1116--1148] where $\mu,\nu\in\R$ and the
boundary condition is not of Feller's type when \mbox{$\nu<0$}. We show
that when $f$ belongs to $C_b^1$ with $f(\infty)=0$ then the
following probabilistic representation of the solution is valid:
\[
u(t,x) = \EE_x \bigl[f(X_t) \bigr] -
\EE_x \biggl[ f'(X_t) \int
_0^{\ell_t^0(X)} e^{-2(\nu-\mu)s} \,ds \biggr],
\]
where $X$ is a reflecting Brownian motion with drift $\mu$ and
$\ell^0(X)$ is the local time of $X$ at $0$. The solution can be
interpreted in terms of $X$ and its creation in $0$ at rate
proportional to $\ell^0(X)$. Invoking the law of $(X_t,\ell_t^0(X))$,
this also yields a closed integral formula for $u$ expressed in
terms of $\mu$, $\nu$ and $f$.
\end{abstract}

% KEYWORDS

\end{frontmatter}

%s1 #&#
\section{Introduction}
%%%%%%%%%%%%%%%%%%%%%%

In this paper, we consider the initial boundary value problem
%
%e1.1 #&#
%e1.2 #&#
%e1.3 #&#
\begin{eqnarray}
\label{1.1} u_t &= &\mu u_x + \tfrac{1}{2}
u_{xx} \qquad(t>0, x \ge0),
\\
\label{1.2} u(0,x) &=& f(x) \qquad(x \ge0),
\\
\label{1.3} u_t(t,0)& =& \nu u_x(t,0)\qquad (t>0)
\end{eqnarray}
of Stroock and Williams \cite{SW-1} (see also \cite{SW-2,WA,PS-1,PS-2})
where $\mu,\nu\in\R$ and the
boundary condition is not of Feller's type when $\nu<0$ (cf. \cite{Fe-1,Fe-2,Fe-3}). If $\nu>0$, then it is known
that the solution to \eqref{1.1}--\eqref{1.3} with $f \in
C_b([0,\infty))$ can be represented as
%
%e1.4 #&#
\begin{equation}
\label{1.4} u(t,x) = \EE_x \bigl[ f(\tilde X_t)
\bigr],
\end{equation}
where $\tilde X$ starts at $x$ under $\PP_{ x}$, behaves like
Brownian motion with drift $\mu$ when in $(0,\infty)$, and exhibits
a sticky boundary behaviour at $0$. The process $\tilde X$ can be
constructed by a familiar time change of the reflecting Brownian
motion $X$ with drift $\mu$ (the inverse of the running time plus
the local time of $X$ at $0$ divided by $\nu$) forcing it to spend
more time at $0$ (cf. \cite{IM-1}, page~186). If $\nu=0$, then
\eqref{1.4} remains valid with $\tilde X$ being absorbed at $0$
(corresponding to the limiting case of infinite stickiness). If
$\nu<0$, then Feller's semigroup\vspace*{1pt} approach (cf. \cite{Hi,Yo,Fe-1,Fe-2,Fe-3}) is no longer applicable since
the speed measure of $\tilde X$ cannot be negative. Stroock and
Williams~\cite{SW-1} show that the minimum principle breaks down in
this case (nonnegative $f$ can produce negative $u$) so that the
solution to \eqref{1.1}--\eqref{1.3} cannot be represented by
\eqref{1.4} where $\tilde X$ is a strong Markov process which
behaves like Brownian motion with drift $\mu$ when in $(0,\infty)$
(for connections with Feller's Brownian motions see \cite{IM-2}, Section~5.7).

Motivated by this peculiarity, Stroock and Williams \cite{SW-1} show
that the solution to \eqref{1.1}--\eqref{1.3} is still generated by a
semi-group of operators when $\nu<0$ and they characterise
nonnegative solutions by means of the Riccati equation. This leads
to subspaces of functions $f$ for which \eqref{1.4} remains valid
with the same time-changed Brownian motion $X$ with drift $\mu$ that
now jumps into $(0,\infty)$ or possibly to a coffin state just
before hitting $0$. This representation of the solution is
applicable when $f(0) = \int_0^\infty f(y) g(y) \,dy$ where $g$ is
the minimal nonnegative solution to the Riccati equation. For more
details and further fascinating developments along these lines, see
\cite{SW-1,SW-2,WA,PS-1,PS-2}.

Inspired by these insights, in this paper we develop an entirely
different approach to solving \eqref{1.1}--\eqref{1.3}
probabilistically that applies to smooth initial data $f$ vanishing
at $\infty$ with no further requirement on its shape. First,
exploiting higher degrees of smoothness of the solution $u$ in the
interior of the domain (which is a well-known fact from the theory
of parabolic PDEs), we reduce the \emph{sticky} boundary behaviour at
$0$ to (i) a \emph{reflecting} boundary behaviour when $\nu= \mu$
and (ii) an \emph{elastic} boundary behaviour when $\nu\ne\mu$.
Second, writing down the probabilistic representations of the
solutions to the resulting initial boundary value problems expressed
in terms of the reflecting Brownian motion with drift $\mu$ and its
local time at~$0$, choosing joint realisations of these processes
where the initial point is given explicitly so that the needed
algebraic manipulations are possible (making use of the extended
L\'evy's distributional theorem), we find that the following
probabilistic representation of the solution is valid:
%
%e1.5 #&#
\begin{equation}
\label{1.5} u(t,x) = \EE_x \bigl[ F \bigl(X_t,
\ell_t^0(X) \bigr) \bigr],
\end{equation}
where $X$ is a reflecting Brownian motion with drift $\mu$ starting
at $x$ under $\PP_{ x}$, and $\ell^0(X)$ is the local time of $X$
at $0$. The function $F$ is explicitly given by
%
%e1.6 #&#
\begin{equation}
\label{1.6} F(x,\ell) = f(x) - f'(x) \int_0^\ell
e^{-2(\nu-\mu)s} \,ds
\end{equation}
for $x \ge0$ and $\ell\ge0$. The derivation applies
simultaneously to all $\mu$ and $\nu$ with no restriction on the
sign of $\nu$, and the process $X$ (with its local time) plays the
role of a fundamental solution in this context (a building block for
all other solutions).

Since $(X,\ell^0(X))$ is a Markov process, we see that the solution
$u$ is generated by the semi-group of transition operators
$(\PP_t)_{t \ge0}$ acting on $f$ by means of \eqref{1.5} and
\eqref{1.6} (in the reverse order). Moreover, it is clear from
\eqref{1.5} and \eqref{1.6} that the solution can be interpreted in
terms of $X$ and its creation in $0$ at rate proportional to
$\ell^0(X)$. Note that this also holds when $\nu<0$ in which case
the Feller's semi-group approach based on the probabilistic
representation \eqref{1.4} is not applicable. Finally, invoking the
law of $(X_t,\ell_t^0(X))$ we derive a closed integral formula for
$u$ expressed in terms of $\mu$, $\nu$ and $f$. Integrating further
by parts yields a closed formula for $u$ where smoothness of $f$ is
no longer needed.

%s2 #&#
\section{Result and proof}
%%%%%%%%%%%%%%%%%%%%%%%%%%

Consider the initial boundary value problem \eqref{1.1}--\eqref{1.3}
and recall that $C_b^1([0,\infty))$ denotes the family of $C^1$
functions $f$ on $[0,\infty)$ such that $f$ and $f'$ are bounded on
$[0,\infty)$. Recall also that the standard normal density and tail
distribution functions are given by $\varphi(x) = (1/\sqrt{2\pi})
e^{-x^2/2}$ and $\Psi(x) = 1 - \Phi(x) = \int_x^\infty\varphi(y)
\,dy$ for $x \in\R$, respectively. The main result of the paper may be
stated as follows.

%%%%%%%%%%%%%%%%%%%%%%%%%%%%%%%%%%%%%%%%%%%%%%%%%%%%%%%%%%%%%%%%%%%%%%%%%%%%%%%%%%
%%% Figure~1 %%%
%%%%%%%%%%%%%%%%%%%%%%%%%%%%%%%%%%%%%%%%%%%%%%%%%%%%%%%%%%%%%%%%%%%%%%%%%%%%%%%%%%

%f1 #&#
\begin{figure}

\includegraphics{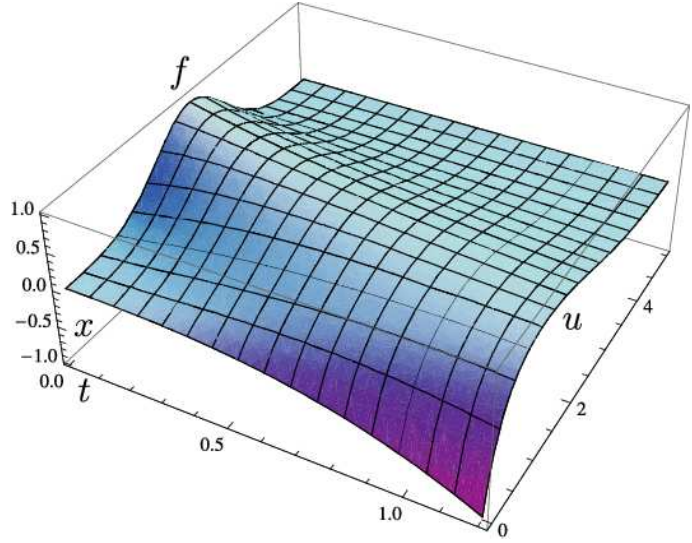}

\caption{The solution $u$ to the initial boundary value
problem (\protect\ref{1.1})--(\protect\ref{1.3}) when $\mu=1$, $\nu=-1/2$ and
$f(x)=\exp(-(x-5/2)^2)$ for $x \ge0$. Note that $u$ takes negative
values even though $f$ is positive so that the classic semi-group
representation \protect\eqref{1.4} of $u$ is not possible in this
case. The
probabilistic representation \protect\eqref{2.1} is valid and this also
yields the integral representation \protect\eqref{2.2}. The solution
can be
interpreted in terms of a reflecting Brownian motion $X$ with drift
$\mu$ and its creation in $0$ at rate proportional to $\ell
^0(X)$.}\label{fig1}
\end{figure}

%%%%%%%%%%%%%%%%%%%%%%%%%%%%%%%%%%%%%%%%%%%%%%%%%%%%%%%%%%%%%%%%%%%%%%%%%%%%%%%%%%

\begin{theorem}\label{th1}
\textup{(i)} If $f \in C_b^1([0,\infty))$ with
$f(\infty)=0$, then there exists a unique solution $u$ to
\eqref{1.1}--\eqref{1.3} satisfying $u \in C^\infty((0,\infty)
\times[0,\infty))$ with $u,u_x \in C_b([0,T] \times
[0,\infty))$ for $T>0$ and $u(t,\infty)=0$ for $t>0$.\vspace*{-6pt}
\begin{longlist}[(iii)]
\item[(ii)] The solution $u$ admits the following probabilistic
representation:
%
%e2.1 #&#
\begin{equation}
\label{2.1} u(t,x) = \EE_x \bigl[f(X_t) \bigr] -
\EE_x \biggl[ f'(X_t) \int
_0^{\ell_t^0(X)} e^{-2(\nu-\mu)s} \,ds \biggr],
\end{equation}
where $X$ is a reflecting Brownian motion with drift $\mu$ starting
at $x$ under $\PP_{ x}$, and $\ell^0(X)$ is the local time of $X$
at $0$ (see Figure \ref{fig1} below).

\item[(iii)] The solution $u$ admits the following integral representation:
%
%e2.2 #&#
\begin{equation}
\label{2.2} u(t,x) = \int_0^\infty f(y) G(t;x,y)
\,dy - \int_0^\infty f'(y) H(t;x,y) \,dy
\end{equation}
where the kernels $G$ and $H$ are given by
%
%e2.3 #&#
%e2.4 #&#
\begin{eqnarray}
\label{2.3} G(t;x,y) &=& \frac{1}{\sqrt{t}} \biggl[ e^{2\mu y} \varphi
\biggl(\frac{x
+ y + \mu t}{\sqrt{t}} \biggr) + \varphi \biggl(\frac{x - y + \mu
t}{\sqrt {t}} \biggr)
\nonumber
\\[-8pt]
\\[-8pt]
\nonumber
&&\hspace*{76pt}{} - 2\mu
e^{2\mu y} \Psi \biggl(\frac{x + y + \mu t}{\sqrt{t}} \biggr) \biggr],
\\
\label{2.4} H(t;x,y) &=& \frac{e^{2\mu y}}{\nu - \mu} \biggl[ (2\nu - \mu)
e^{2(
\nu-\mu)(x + y + \nu t)} \Psi \biggl(\frac{x + y + (2 \nu - \mu)t}{
\sqrt{t}} \biggr)\nonumber
\\
&&\hspace*{155pt}{}- \mu\Psi \biggl(\frac{x + y + \mu t}{\sqrt{t}} \biggr)
 \biggr]\nonumber \\
 \\[-20pt]
 \eqntext{\mbox{if } \nu\ne
\mu}
\\
\nonumber
&=& 2 e^{2\mu y} \biggl[ \bigl(1 + \mu(x + y + \mu t) \bigr) \Psi
\biggl(\frac{x + y + \mu t}{\sqrt{t}} \biggr) \\
&&\hspace*{90pt}{}- \mu\sqrt{t} \varphi \biggl(\frac{x + y + \mu t}{\sqrt{t}}
\biggr) \biggr] \qquad\mbox{if } \nu =\mu\nonumber
\end{eqnarray}
for $t>0$ and $x,y \ge0$.
\end{longlist}
\end{theorem}

\begin{pf} Let $f \in C_b^1([0,\infty))$ with $f(\infty)=0$ be
given and fixed. We first show that any solution $u$ to
\eqref{1.1}--\eqref{1.3} satisfying $u \in C^\infty((0,\infty)
\times[0,\infty))$ with $u,u_x \in C_b([0,T] \times
[0,\infty))$ for $T>0$ and $u(t,\infty)=0$ for $t>0$ admits the
probabilistic representation \eqref{2.1}.

\begin{longlist}[1.]
\item[1.] Setting $v=u_x$ and differentiating both sides in \eqref{1.1}
with respect to $x$ we see that $v$ solves the same equation
%
%e2.5 #&#
\begin{equation}
\label{2.5} v_t = \mu v_x + \tfrac{1}{2}
v_{xx} \qquad(t>0, x \ge0).
\end{equation}
Moreover, differentiating both sides in \eqref{1.2} with respect to
$x$ we find that
%
%e2.6 #&#
\begin{equation}
\label{2.6} v(0,x) = f'(x)\qquad (x \ge0).
\end{equation}
Finally, combining \eqref{1.3} with \eqref{1.1} we see that
\eqref{1.3} reads as follows:
%
%e2.7 #&#
\begin{equation}
\label{2.7} v_x(t,0) = \lambda v(t,0) \qquad(t>0),
\end{equation}
where we set $\lambda= 2(\nu - \mu)$. In this way, we have obtained
the initial boundary value problem \eqref{2.5}--\eqref{2.7} for $v$.
Note that the boundary condition \eqref{2.7} corresponds to
(i)~a\vadjust{\goodbreak}
reflecting boundary behaviour when $\lambda=0$ and (ii) an elastic
boundary behaviour when $\lambda\ne0$. Setting
%
%e2.8 #&#
\begin{equation}
\label{2.8} B_t^{-\mu} = B_t - \mu t\quad \mbox{and}\quad
S_t^{-\mu} = \sup_{0 \le s \le t}
B_s^{-\mu}
\end{equation}
for $t \ge0$ where $B$ is a standard Brownian motion, and denoting
by $R^{\mu,x}$ a reflecting Brownian motion with drift $\mu$
starting at $x$ in $[0,\infty)$, it is known that the classic
L\'evy's distributional theorem (see \cite{RY}, page~240) extends as
follows:
%
%e2.9 #&#
\begin{equation}
\label{2.9} \bigl(x \vee S^{-\mu} - B^{-\mu},x \vee
S^{-\mu} - x \bigr) \buildrel \mathrm{law} \over= \bigl(R^{\mu,x},
\ell^0 \bigl(R^{\mu,x} \bigr) \bigr),
\end{equation}
where $\ell^0(R^{\mu,x})$ is the local time of $R^{\mu,x}$ at $0$
(for a formal verification based on Skorokhod's lemma see the proof
of Theorem~3.1 in \cite{Pe}). Identifying
%
%e2.10 #&#
\begin{equation}
\label{2.10} X_t^x:= x \vee S_t^{-\mu}
- B_t^{-\mu}\quad \mbox{and}\quad \ell_t^0
\bigl(X^x \bigr) = x \vee S_t^{-\mu} - x
\end{equation}
in accordance with \eqref{2.9} above, we claim (cf. \cite{IM-1}, pages~183--184) that the solution $v$ to the problem
\eqref{2.5}--\eqref{2.7} admits the probabilistic representation
%
%e2.11 #&#
\begin{equation}
\label{2.11} v(t,x) = \EE \bigl[ e^{-\lambda\ell_t^0(X^x)} f'
\bigl(X_t^x \bigr) \bigr]
\end{equation}
for $t \ge0$ and $x \ge0$ (for multi-dimensional extensions see
\cite{BCS}, Section~2).

\item[2.] To verify \eqref{2.11}, we can make use of standard arguments by
letting time run backward and applying It\^o's formula to $v$
composed with $(t - s,X_s^x)$ and multiplied by $e^{-\lambda
\ell_s^0(X^x)}$ for $s \in[0,t)$ where $t>0$ and $x \ge0$ are
given and fixed. This yields
%
%e2.12 #&#
\begin{eqnarray}
\label{2.12} &&e^{-\lambda\ell_s^0(X^x)} v \bigl(t - s,X_s^x
\bigr)\nonumber\\
&&\qquad= v(t,x) + \int_0^s (-\lambda)
e^{-\lambda\ell_r^0(X^x)} v \bigl(t - r,X_r^x \bigr) \,d
\ell_r^0 \bigl(X^x \bigr)\nonumber
\\
\nonumber
&&\qquad\quad{}+ \int_0^s e^{-\lambda\ell_r^0(X^x)}
(-v_t) \bigl(t - r,X_r^x \bigr) \,dr
\\
\nonumber
&&\qquad\quad{}+ \int_0^s e^{-\lambda\ell_r^0(X^x)}
v_x \bigl(t - r,X_r^x \bigr) \,d \bigl(x \vee
S_r^{-\mu} - B_r^{-\mu} \bigr)
\\
&&\qquad\quad{}+ \frac{1}{2} \int_0^s
e^{-\lambda\ell_r^0(X^x)} v_{xx} \bigl(t - r,X_r^x
\bigr) \,d \bigl\langle X^x, X^x \bigr\rangle_r
\\
\nonumber
&&\qquad= v(t,x) + \int_0^s
e^{-\lambda\ell_r^0(X^x)} (-\lambda v + v_x) \bigl(t - r,X_r^x
\bigr) \,d \bigl(x \vee S_r^{-\mu} \bigr)
\\
\nonumber
&&\qquad\quad{}+ \int_0^s e^{-\lambda\ell_r^0(X^x)}
\biggl(-v_t + \mu v_x + \frac{1} {2}
v_{xx} \biggr) \bigl(t - r,X_r^x \bigr) \,dr
\\
\nonumber
&&\qquad\quad{}- \int_0^s e^{-\lambda
\ell_r^0(X^x)}
v_x \bigl(t - r,X_r^x \bigr)
\,dB_r
\\
\nonumber
&&\qquad= v(t,x) - \int_0^s
e^{-\lambda\ell_r^0(X^x)} v_x \bigl(t - r,X_r^x
\bigr) \,dB_r
\end{eqnarray}
since $d(x \vee S_r^{-\mu})$ is zero off the set of all $r$ at
which $X_r^x \ne0$, while $(-\lambda v + v_x)(t - r,X_r^x) = 0$
for $X_r^x=0$ by \eqref{2.7} above, so that the integral with
respect to $d(x \vee S_r^{-\mu})$ is equal to zero. Note also
that $d \langle X^x, X^x \rangle_r = dr$ since $r \mapsto x \vee
S_r^{-\mu}$ is increasing, and thus of bounded variation while in
the final equality we also use~\eqref{2.5}. From \eqref{2.12}, we see
that
%
%e2.13 #&#
\begin{equation}
\label{2.13} v(t,x) = e^{-\lambda\ell_s^0(X^x)} v \bigl(t - s,X_s^x
\bigr) + M_s,
\end{equation}
where $M_s = \int_0^s e^{-\lambda\ell_r^0(X^x)} v_x(t - r,X_r^x)
\,dB_r$ is a continuous local martingale for $s \in[0,t)$. Choose a
localisation sequence of stopping times $(\sigma_n)_{n \ge1}$ for
$M$ (meaning that $M$ stopped at $\sigma_n$ is a martingale for each
$n \ge1$ and $\sigma_n \uparrow\infty$ as $n \rightarrow\infty$),
take any sequence $s_n \uparrow t$ as $n \rightarrow\infty$, and
set $\tau_n:= \sigma_n \wedge s_n$ for $n \ge1$. Then the optional
sampling theorem yields
%
%e2.14 #&#
\begin{eqnarray}
\label{2.14} v(t,x) &=& \EE \bigl[ e^{-\lambda\ell_{\tau_n}^0(X^x)} v \bigl(t -
\tau_n, X_{\tau_n}^x \bigr) \bigr] + \EE
M_{\tau_n}
\nonumber\\
&=& \EE \bigl[ e^{- \lambda\ell_{\tau_n}^0(X^x)}v \bigl(t - \tau_n,X_{\tau_n}^x
\bigr) \bigr] \rightarrow\EE \bigl[ e^{-\lambda\ell_t^0(X^x)} v \bigl(0,X_t^x
\bigr) \bigr] \\
&= &\EE \bigl[ e^{-\lambda\ell_t^0(X^x)}f' \bigl(X_t^x
\bigr) \bigr]\nonumber
\end{eqnarray}
as $n \rightarrow\infty$ by the dominated convergence theorem and
\eqref{2.6} above where we use that $v \in C_b([0,T] \times
[0,\infty))$ for $T \ge t$ and $\EE e^{\vert\lambda\vert
\ell_t^0(X^x)} < \infty$ for $t>0$ in view of~\eqref{2.10} above.
This establishes \eqref{2.11} as claimed.

\item[3.] Recalling that $v=u_x$ and $u(t,\infty)=0$ we find using
\eqref{2.10} and \eqref{2.11} that
%
%e2.15 #&#
\begin{eqnarray}
\label{2.15} u(t,x) &=& - \int_x^\infty
u_x(t,y) \,dy + u(t,\infty)\nonumber
\\
\nonumber
&=& - \int_x^\infty v(t,y) \,dy \\
&=& - \int
_x^\infty\EE \bigl[ e^{-\lambda
(y \vee S_t^{-\mu}-y)} f'
\bigl(y \vee S_t^{- \mu} - B_t^{-\mu}
\bigr) \bigr] \,dy\nonumber
\\
\nonumber
&=& - \int_x^\infty\EE \bigl[
f' \bigl(y - B_t^{-\mu} \bigr) I
\bigl(S_t^{-\mu} \le y \bigr) \\
&&\hspace*{40pt}{}+ e^{-\lambda(S_t^{-\mu}-y)}
f' \bigl(S_t^{-\mu} - B_t^{-\mu}
\bigr) I \bigl(S_t^{-\mu} > y \bigr) \bigr] \,dy
\nonumber
\\[-8pt]
\\[-8pt]
\nonumber
&=& - \EE \biggl[\int_{x \vee S_t^{-\mu}}^\infty
f' \bigl(y - B_t^{-\mu} \bigr) \,dy \biggr]
\\
&&{}- \EE
\biggl[ \int_x^{x \vee S_t^{-\mu}} e^{-\lambda(S_t^{-\mu}
-y)}f'
\bigl(S_t^{-\mu} - B_t^{-\mu} \bigr) \,dy
\biggr]\nonumber
\\
&=& - \EE \biggl[\int_{x \vee S_t^{-\mu} - B_t^{-\mu}}^\infty
f'(z) \,dz \biggr] \nonumber\\
&&{}- \EE \biggl[ f' \bigl(x \vee
S_t^{-\mu} - B_t^{-\mu} \bigr) \int
_x^{x \vee
S_t^{-\mu}} e^{-\lambda(x \vee S_t^{-\mu}-y)} \,dy \biggr]
\nonumber\\
\nonumber
&=& \EE \bigl[f \bigl(x \vee S_t^{-\mu} -
B_t^{-\mu} \bigr) \bigr] - \EE \biggl[ f'
\bigl(x \vee S_t^{-\mu} - B_t^{-\mu}
\bigr) \int_0^{x \vee S_t^{-\mu}
-x} e^{-\lambda s} \,ds \biggr]
\end{eqnarray}
for $t \ge0$ and $x \ge0$, where in the second last equality we
use that $S_t^{-\mu} = x \vee S_t^{-\mu}$ since otherwise the
integral from $x$ to $x \vee S_t^{-\mu}$ equals zero, and in the
last equality we use that $f(\infty)=0$. Making use of \eqref{2.9}
in \eqref{2.15} establishes the probabilistic representation
\eqref{2.1} as claimed in the beginning of the proof.

%%% End Figure~4.\
\item[4.] Focusing on \eqref{2.1} and recalling \eqref{2.10}, we see that
an explicit calculation of the right-hand side in \eqref{2.1} is
possible since the probability density function $g$ of
$(B_t^{-\mu},S_t^{-\mu})$ is known and can be readily derived from
the known probability density function of $(B_t,S_t)$ when $\mu$ is
zero (see, e.g., \cite{IM-2}, page~27 or \cite{RY}, page~110) using a
standard change-of-measure argument. This yields the following
closed form expression:
%
%e2.16 #&#
\begin{equation}
\label{2.16} g(t;b,s) = \sqrt{\frac{2}{\pi}} \frac{1}{t^{3/2}} (2s - b)
\exp \biggl[ - \frac{(2s - b)^2}{2t} - \mu \biggl(b + \frac{\mu
t}{2} \biggr)
\biggr]
\end{equation}
for $t>0$ and $b \le s$ with $s \ge0$. It follows that the
functions on the right-hand side of \eqref{2.1} can be given the
following integral representations:
%
%e2.17 #&#
%e2.18 #&#
\begin{eqnarray}
\label{2.17} u^1(t,x) &:= &\EE_x \bigl[f(X_t)
\bigr] = \EE \bigl[ f \bigl(x \vee S_t^{-\mu} -
B_t^{-\mu} \bigr) \bigr]
\nonumber
\\[-8pt]
\\[-8pt]
\nonumber
&=& \int_0^
\infty\int
_{-\infty}^s f(x \vee s - b) g(t;b,s) \,db \,ds,
\\
\label{2.18} u^2(t,x) &:=& \EE_x \biggl[
f'(X_t) \int_0^{\ell_t^0(X)}
e^{-\lambda r} \,dr \biggr]
\nonumber\\
&=& \EE_x \biggl[ f' \bigl(x \vee
S_t^{-\mu} - B_t^{-\mu} \bigr) \int
_0^{x \vee S_t^{-\mu} - x} e^{-\lambda r} \,dr \biggr]
\\
\nonumber
&=& \int_0^\infty\int
_{-\infty}^s \biggl( f'(x \vee s - b) \int
_0^{x \vee s - x} e^{-\lambda r} \,dr \biggr) g(t;b,s) \,db \,ds
\end{eqnarray}
for $t>0$ and $x \ge0$ where $\lambda= 2(\nu - \mu)$. A lengthy
elementary calculation then shows that
%
%e2.19 #&#
%e2.20 #&#
\begin{eqnarray}
\label{2.19} u^1(t,x) &=& \int_0^\infty
f(y) G(t;x,y) \,dy,
\\
\label{2.20} u^2(t,x) &=& \int_0^\infty
f'(y) H(t;x,y) \,dy
\end{eqnarray}
for $t>0$ and $x \ge0$ where $G$ and $H$ are given in \eqref{2.3}
and \eqref{2.4} above. Noting that
%
%e2.21 #&#
\begin{equation}
\label{2.21} u(t,x) = u^1(t,x) - u^2(t,x)
\end{equation}
we see that this establishes the integral representation \eqref{2.2}
as claimed.

\item[5.] A direct analysis of the integral representations \eqref{2.19}
and \eqref{2.20} with $G$ and $H$ from \eqref{2.3} and \eqref{2.4}
then shows that $u$ from \eqref{2.21} belongs to both
$C^\infty((0,\infty) \times[0,\infty))$ and $C_b([0,T]
\times[0,\infty))$ for $T>0$ and $u(t,\infty)=0$ for $t>0$. A
similar analysis also shows that both $u_x^1$ and $u_x^2$ belong to
$C_b(([0,T] \times[0,\infty)) \setminus\{(0,0)\})$ for
$T>0$. Moreover, it can be directly verified that (i) $u_x^1(t,x)
\rightarrow f'(x)$ as $t \downarrow0$ for all $x>0$ but
$u_x^1(t,0)=0$ for all $t>0$ so that $u_x^1$ is not continuous at
$(0,0)$ unless $f'(0)=0$; and (ii) $u_x^2(t,x) \rightarrow0$ as $t
\downarrow0$ for all $x>0$ but $u_x^2(t,0) \rightarrow-f'(0)$ as
$t \downarrow0$ so that $u_x^2$ is not continuous at $(0,0)$ either
unless $f'(0)=0$. Despite the possibility that both $u_x^1$ and
$u_x^2$ are discontinuous at $(0,0)$, it turns out that when acting
in cohort to form $u_x = u_x^1 - u_x^2$ the resulting function
$u_x$ is continuous at $(0,0)$ so that $u_x$ belongs to $C_b([0,T]
\times[0,\infty))$ for $T>0$. It follows therefore from the
construction and these arguments that the function $u$ defined by
\eqref{2.2} with $G$ and $H$ from \eqref{2.3} and \eqref{2.4} solves
the initial boundary problem \eqref{1.1}--\eqref{1.3} and satisfies
$u \in C^\infty((0,\infty) \times[0,\infty))$ with $u,u_x \in
C_b([0,T] \times[0,\infty))$ for $T>0$ and $u(t,\infty)=0$ for
$t>0$. Placing then any such $u$ at the beginning of the proof and
repeating the same arguments as above, we can conclude that $u$ admits
the probabilistic representation \eqref{2.1}. These arguments
therefore establish both the existence and uniqueness of the
solution $u$ to the initial boundary problem \eqref{1.1}--\eqref{1.3}
satisfying the specified conditions and the proof is complete.\quad\qed
\end{longlist}
\noqed\end{pf}

\begin{remark}[(Nonsmooth initial data)] The integral
representation \eqref{2.2} requires that $f$ is differentiable.
Integrating by parts we find that
%
%e2.22 #&#
\begin{equation}
\label{2.29} \qquad\int_0^\infty f'(y)
H(t;x,y) \,dy = -f(0) H(t;x,0) - \int_0^\infty f(y)
H_y(t;x,y) \,dy.
\end{equation}
Inserting this back into \eqref{2.2}, we find that $u$ admits the
following integral representation:
%
%e2.23 #&#
\begin{equation}
\label{2.30} u(t,x) = \int_0^\infty f(y) (G +
H_y) (t;x,y) \,dy + f(0) H(t;x,0),
\end{equation}
where the first function is given by
%
%e2.24 #&#
\begin{eqnarray}
\label{2.31}&& (G + H_y) (t;x,y)\nonumber\\
&&\qquad= \frac{1}{\sqrt{t}} \biggl[
\varphi \biggl( \frac{x - y + \mu t}{\sqrt{t}} \biggr) - e^{2\mu y} \varphi \biggl(
\frac{x + y + \mu t}{\sqrt{t}} \biggr) \biggr]\nonumber
\\
&&\qquad\quad{}- \frac{2 \nu e^{2\mu y}} {\nu - \mu} \biggl[ \mu\Psi \biggl(\frac{x + y + \mu t} {\sqrt{t}}
\biggr)\\
&&\hspace*{85pt}{} + (\mu - 2\nu) e^{2(\nu-\mu)(x+y+\nu t)}
 \Psi \biggl(\frac{x + y + (2\nu - \mu) t}{\sqrt{t}} \biggr) \biggr]
 \nonumber\\
 \eqntext{\mbox{if } \nu
\ne \mu}
\\
\nonumber
&&\qquad= \frac{1}{
\sqrt{t}} \varphi \biggl(\frac{x - y + \mu t}{\sqrt{t}} \biggr)\\
&&\qquad\quad{} -
\frac{e^{2\mu y}}{\sqrt{t}} \biggl[ \bigl(1 + 4 \mu^2 t \bigr) \varphi
\biggl(
\frac{x + y + \mu t}{\sqrt{t}} \biggr)\nonumber
\\
&&\hspace*{73pt}{}- 4 \mu \bigl(1 + \mu(x + y) + \mu^2 t \bigr) \sqrt{t}
\Psi \biggl(\frac{x + y + \mu t}{\sqrt{t}} \biggr) \biggr]\nonumber\\
\eqntext{\mbox{if } \nu= \mu}
\end{eqnarray}
and the second function is given by
%
%e2.25 #&#
\begin{eqnarray}
\label{2.32} H(t;x,0)&= &\frac{1}{\nu - \mu} \biggl[ (2\nu - \mu)
e^{2(
\nu-\mu)(x + \nu t)} \Psi \biggl(\frac{x + (2 \nu - \mu)t}{
\sqrt{t}} \biggr)\nonumber\\
&&\hspace*{145pt}{} - \mu\Psi \biggl(
\frac{x + \mu t}{\sqrt{t}} \biggr) \biggr] \qquad\mbox{if } \nu\ne\mu
\\
\nonumber
&= & 2 \biggl[ \bigl(1 + \mu(x + \mu t) \bigr) \Psi \biggl(
\frac{x +
\mu t}{\sqrt{t}} \biggr) - \mu\sqrt{t} \varphi \biggl(\frac{x + \mu
t}{\sqrt{t}} \biggr)
\biggr] \qquad\mbox{if } \nu=\mu\nonumber
\end{eqnarray}
for $t>0$ and $x,y \ge0$. Note that smoothness of $f$ is no longer
needed in the integral representation \eqref{2.30} and this formula
for $u$ can be used when $f \in C_b([0,\infty))$ for instance.
\end{remark}

%%%%%%%%%%%%%%%%%%%%%%%%%%%%%%%%%%%%%%%%%%%%%%%%%%%%%%%%%%%%%%%%%%%%%%%%%%%%%%%%%%
%%% References %%%
%%%%%%%%%%%%%%%%%%%%%%%%%%%%%%%%%%%%%%%%%%%%%%%%%%%%%%%%%%%%%%%%%%%%%%%%%%%%%%%%%%

% imsref loaded by akundreckaite, 2014-01-31 08:30:56

%%%%%%%%%%%%%%%%%%%%%%%%%%%%%%%%%%%%%%%%%%%%%%%%%%%%%%%%%%%%%%%%%%%%%%%%%%%%%%%%%%
%%% Affiliations %%%
%%%%%%%%%%%%%%%%%%%%%%%%%%%%%%%%%%%%%%%%%%%%%%%%%%%%%%%%%%%%%%%%%%%%%%%%%%%%%%%%%%


%suskaldyti doi

\printaddresses

\end{document}